\documentclass[12pt,a4paper]{article}
\usepackage{latexsym,amsmath,amssymb,amscd}
\newtheorem{theorem}{Theorem}

\author{J.~Brodzki \and R.~J.~Plymen}
\title{Geometry of the smooth dual of $GL(n)$}
\date{}
\begin{document}
\maketitle
%\markboth{J.~Brodzki, R.~J.~Plymen}{Geometry of the smooth dual of $GL(n)$ }

\begin{abstract}
Let $\mathcal{A}(n)$ be the smooth dual of the $p$-adic group $G=GL(n)$. 
We create on $\mathcal{A}(n)$ the structure of a complex algebraic variety. 
There is a morphism of $\mathcal{A}(n)$ onto the Bernstein variety 
$\Omega G$ which is injective on each component of $\mathcal{A}(n)$. 
The tempered dual of $G$
is a deformation retract of $\mathcal{A}(n)$. 
The periodic cyclic homology of
the Hecke algebra of $G$  is isomorphic to the periodised de Rham cohomology 
 supported on finitely many components of $\mathcal{A}(n)$.  
\end{abstract}

\begin{center}
\bf Introduction
\end{center}

Let $F$ be a finite extension of $\mathbb{Q}_p$  and let $G=GL(n) = GL(n,F)$. 
A representation of $G$ on a complex vector space $V$ is {\it smooth}
if the stabilizer of each vector in  $V$ is an open subgroup of $G$. 
The set of equivalence classes of irreducible smooth representations 
of $G$ is the {\it smooth dual} $\mathcal{A}(n)$ of $G$. The Bernstein variety $\Omega
G$ is the set of $G$-conjugacy classes of pairs $(M, \sigma)$ where 
$M$ is a Levi subgroup of $G$ and $\sigma$ is an irreducible supercuspidal 
representation of $M$. Each irreducible smooth representation of $G$ 
is a subquotient of an induced representation $i_{GM}\sigma$. The 
pair $(M, \sigma)$ is unique up to conjugacy. This creates a finite-to-one 
map, the infinitesimal character,  
 from $ \mathcal{A}(n)$ onto $\Omega G$.  

The set $\Omega G$ has the structure of a complex algebraic variety. 
This structure comes ultimately from twisting a supercuspidal 
representation $\sigma$ of $M$ by  unramified quasicharacters of $M$. 
The set $\Psi(M)$ of all unramified quasicharacters
of $M$  has the structure of a
complex torus.

The local Langlands correspondence \cite{HT}\cite{H} is a natural bijection 
from $\mathcal{G}(n)$ to $\mathcal{A}(n)$ where $\mathcal{G}(n)$ 
is the set of equivalence classes of admissible $n$-dimensional 
complex representations of the Weil-Deligne group $W'_F$. 

Each representation of the Weil group $W_F$ may itself be twisted 
by an unramified quasicharacter of $W_F$. In this way, by adapting 
Bernstein's construction, we create the structure of complex algebraic
variety on $\mathcal{G}(n)$. We then transport this structure 
to $\mathcal{A}(n)$ by the local Langlands correspondence. So $\mathcal{A}(n)$ 
acquires the structure of complex algebraic variety with infinitely 
many components. We obtain the following results. 

\begin{theorem}
The smooth dual $\mathcal{A}(n)$ has a natural structure of complex 
algebraic variety. There is a canonical morphism of 
$\mathcal{A}(n)$ onto the Bernstein variety $\Omega G$ which is injective
on each component of $\mathcal{A}(n)$. 
\end{theorem}
The Langlands parameters effect a stratification of the smooth dual $\mathcal{A}(n)$. 
\begin{theorem}
The tempered dual of $G$  is a deformation retract of $\mathcal{A}(n)$ . 
\end{theorem}
We view the deformation retraction as the geometric counterpart of the assembly 
map in the Baum-Connes conjecture for $GL(n)$, see \cite{BHP1}\cite{L}.

Let $\mathcal{H}(G)$ be the Hecke algebra of compactly supported, uniformly locally 
constant complex-valued functions on $G$. 
\begin{theorem}
The periodic cyclic homology of $\mathcal{H}(G)$ is isomorphic to the 
periodised de Rham  cohomology  supported on
finitely many  components of the smooth dual of $G$. 
\end{theorem}
From the point of view of noncommutative geometry, the variety 
$\mathcal{A}(n)$ is an {\it affine scheme} underlying the non-unital noncommutative
ring $\mathcal{H}(G)$. 

We find it interesting that the complex algebraic variety structure 
on $\mathcal{A}(n)$ has become visible only after employing the local Langlands 
correspondence, in other words only after paying attention to the 
number theory of the local field $F$.

\mbox{}

\noindent{\bf 1. The Complex Algebraic Variety $\mathcal{G}(n)$}

The Weil group $W_F$ fits into a short exact sequence
$
0 \rightarrow I_F \rightarrow W_F \stackrel{d}{\rightarrow} \mathbb{Z}
\rightarrow 0, 
$
where $I_F$ is the inertia group of $F$.
The local Langlands correspondence \cite{HT}\cite{H} asserts
that for all $n \geq 1$
there exists a natural bijection
\[\pi_F :  \mathcal{G}^0(n) \longrightarrow \mathcal{A}^0(n)\]
where $\mathcal{G}^0(n)$ is the set of equivalence classes of irreducible
continuous $n$-dimensional complex representations of
the Weil group $W_F$, and $\mathcal{A}^0(n)$ is the set of
equivalence classes of irreducible supercuspidal representations of
$GL(n)$.

We identify the elements of the set $\mathcal{G}^0(1)$, the quasicharacters of
$W_F$, with quasicharacters of $F^{\times}$ via the reciprocity isomorphism
$r_F : W^{ab}_F \cong F^{\times}$.   The local Langlands correspondence
is compatible with twisting by quasicharacters, see \cite[1.2]{H}.

   A complex representation
of the Weil-Deligne group $W_F'$ is a pair $(\rho,N)$ consisting of a continuous
representation $\rho : W_F \longrightarrow GL(V), dim_{\mathbb{C}}V
= n$, together with a nilpotent endomorphism $N \in End \, V$ such
that $\rho(w)N\rho(w)^{-1} = ||w||N$.   A representation
$\rho' = (\rho,N)$ is {\it admissible} if
$\rho$ is semisimple.   For any $n \geq 1$ the representation
$sp(n)$ is defined by $V = \mathbb{C}^n = \mathbb{C}e_0 + \ldots +
\mathbb{C}e_{n-1}, \rho(w)e_i = ||w||^ie_i$ and $Ne_i = e_{i+1}
(0 \leq i < n-1), Ne_{n-1} = 0$.   

Let $\mathcal{G}(n)$ be the set of equivalence classes of admissible
$n$-dimensional representations of the Weil-Deligne group $W'_F$, and let
$\mathcal{A}(n)$ be the set of equivalence classes of irreducible smooth
representations of $GL(n)$.   To any {\it indecomposable} representation
$\rho \otimes sp(r)$ of $W'_F$ associate the essentially square-integrable
representation $Q(\Delta)$ with $\Delta = \{\pi_F(\rho), |\;\;|\pi_F(\rho), \ldots,
|\;\;|^{r-1}\pi_F(\rho)\}$.   To any {\it admissible} representation
$\rho' = (\rho_1 \otimes sp(r_1)) \oplus \ldots \oplus
(\rho_m \otimes sp(r_m))$
of $W'_F$ associate the Langlands quotient $Q(\Delta_1, \ldots, \Delta_m)$
where $\Delta_j$ is the segment
\[\Delta_j = \{\pi_F(\rho_j), \ldots,|\;\;|^{r_j - 1}\pi_F(\rho_j)\}\]

This creates as in \cite[p. 381]{Ku} a natural bijection

\[\pi_F : \mathcal{G}(n) \longrightarrow \mathcal{A}(n)\]

A quasicharacter $\psi : W_F \longrightarrow \mathbb{C}^{\times}$
is {\it unramified} if $\psi$ is trivial on the inertia group
$I_F$.   In that case we have
$\psi(w) = z^{d(w)}$
with $z \in \mathbb{C}^{\times}$.   The group of unramified quasicharacters
of $W_F$ is denoted $\Psi(W_F)$.   We have $\Psi(W_F) \cong \mathbb{C}^
{\times}$.

Let now
$\rho' = \rho_1 \otimes sp(r_1) \oplus \ldots \oplus \rho_m \otimes
sp(r_m)$ be an admissible representation of $W_F'$.   The set
\[
\{\psi_1\rho_1 \otimes sp(r_1) \oplus \ldots \oplus \psi_m\rho_m
\otimes sp(r_m) :  \psi_1, \ldots, \psi_m \in \Psi(W_F)\}
\]
will be called the {\it orbit} of $\rho'$ under the action of
$(\Psi(W_F))^m = \Psi(W_F) \times \ldots \times \Psi(W_F)$ ( $m$ factors ).
Then $\mathcal{G}(n)$ admits a partition into orbits.   Since
$(\Psi(W_F))^m = (\mathbb{C}^{\times})^m$, a complex torus, each orbit
has the structure of a complex algebraic variety.   In this way, the
set of admissible $n$-dimensional representations of
the Weil-Deligne group acquires the structure of a
complex algebraic variety.   Each component in this variety is the
quotient of a complex torus by a product of symmetric groups.

Let $\Omega G$ be the Bernstein variety of $G$.   Each point in
$\Omega G$ is a conjugacy class of cuspidal pairs $(M,\sigma)$.
A quasicharacter $\psi : M \longrightarrow \mathbb{C}^{\times}$
is {\it unramified} if $\psi$ is trivial on $M^{\circ}$.   The group
of unramified quasicharacters of $M$ is denoted $\Psi(M)$.   We have
$\Psi(M) \cong (\mathbb{C}^{\times})^\ell$ where $\ell$ is the parabolic
rank of $M$.   The group $\Psi(M)$ now creates orbits:
the orbit of $(M,\sigma)$ is $\{(M, \psi \otimes \sigma): \psi \in
\Psi(M)\}$.   Denote this orbit by $D$, and set
$\Omega = D/W(M,D)$, where $W(M)$ is the Weyl group of $M$ and
$W(M,D)$ is the subgroup of $W(M)$ which leaves $D$ globally invariant.
The orbit $D$ has the structure of a complex torus, and so $\Omega$
is a complex algebraic variety.   We view $\Omega$ as
a component in the algebraic variety $\Omega G$.

We recall the {\it extended quotient}.   Let the finite group $\Gamma$
act on the space $X$.   Let $\widehat{X} =
\{(x,\gamma) : \gamma x = x\}$, let $\Gamma$ act on $\widehat{X}$ by
$\gamma_1(x,\gamma) = (\gamma_1x, \gamma_1 \gamma \gamma_1^{-1})$. Then 
$\widehat{X}/\Gamma$ is the  extended quotient of $X$ by $\Gamma$. 
There is a canonical projection $\hat{X}/\Gamma \rightarrow X/\Gamma$.

The Bernstein variety $\Omega G$ is the disjoint union 
of ordinary quotients. We now 
replace the ordinary quotient by the extended quotient to create 
a new variety $\Omega^+ G$. So we have
$$
\Omega G = \bigsqcup D/W(M, D)\; \; \mbox{\rm and}
\; \;   \Omega^+G = \bigsqcup \widehat{D}/W(M, D)
$$

Let $\rho_1, \ldots, \rho_m$ be irreducible complex representations of
$W_F$ with $dim_{\mathbb{C}} \, \rho_j = a_j$ and let
\[
\Delta_j = \{\pi_F(\rho_j), \ldots, |\;\;|^{r_j - 1}\pi_F(\rho_j)\} \; 
\mbox{\rm and}\; \;
R(\Delta_j) = \pi_F(\rho_j) \otimes \ldots \otimes | \;\;|^{r_j - 1}
\pi_F(\rho_j)\]

\setcounter{theorem}{0}
\begin{theorem}  
Let $\pi_F$ be the local Langlands correspondence and let $inf.ch.$ be the infinitesimal 
character of Bernstein \cite[(III.4) p.75]{B}. Then we have 
 a commutative diagram
\[
\begin{CD}
\mathcal{G}(n)      @>{\pi_F}>>     \mathcal{A}(n)\\
@V{\alpha}VV                       @VV{inf.ch}V\\
\Omega^{+}G      @>>{\beta}>     \Omega G
\end{CD}
\]
in which the maps $\alpha$ and $\beta$ are as follows
 \[
\begin{CD}
\rho_1 \otimes sp(r_1) \oplus
\ldots \oplus \rho_m \otimes sp(r_m)  @>{\pi_F}>> Q(\Delta_1, \ldots,\Delta_m)\\
@V{\alpha}VV                                @VV{inf.ch.}V\\
(\prod GL(a_j)^{r_j}, \bigotimes \pi_F(\rho_j)^{\otimes r_j})   @>>{\beta}>
(\prod GL(a_j)^{r_j}, \bigotimes R(\Delta_j))
\end{CD}
\]
The map $\beta : \Omega^{+}G \longrightarrow \Omega G$
is not the canonical projection, but a twisted projection.
The map $\beta$
is a morphism of complex algebraic varieties, injective on each component
of $\Omega^+G$.

 Let $\mathcal{G}_\Omega(n) =
(\beta \circ \alpha)^{-1}\Omega$.  In the Bernstein decomposition
$\mathcal{G}(n) = \bigsqcup \mathcal{G}_\Omega(n)$ we have
\[
\mathcal{G}_\Omega(n) = \Omega^{+}
\]
where we write $\Omega^+ = \widehat{D}/W(M,D)$ when $\Omega = D/W(M,D)$. 
\end{theorem}

Proof.     We have to delve into some combinatorics.
Let $dim_{\mathbb{C}} \rho = a$, let $(M,\sigma) = (GL(a)^r, \pi_F
(\rho)^{\otimes r})$, let $D$ be the orbit of $(M,\sigma)$.
Then $W(M,D) = S_r$ the symmetric group on $r$ letters, and
$\Omega = \mathbb{C}^{\times r}/S_r = Sym^{r} \mathbb{C}^{\times}$, 
the $r$-fold symmetric product of $\mathbb{C}^\times$. 

Write $\Gamma = S_r$.   Let
$
t_1 ,  \ldots ,  t_k 
$
be a partition of $r$ in which $t_1 < \ldots < t_k$ and $t_j$
is repeated $n_j$ times, $1 \leq j \leq k$.   Then we have
$n_1t_1 + \ldots + n_kt_k = r$.  Let $\gamma$ be the permutation
of $r$ letters which corresponds to this partition, so that $\gamma$
is the product of cycles $(1, \ldots, t_1) \ldots (1, \ldots, t_k)$.

Let
\[
\rho'_{\gamma} = \rho \otimes sp(t_1) \oplus \ldots \oplus \rho
\otimes sp(t_k)
\]
The first summand is repeated $n_1$ times, \ldots, the last summand is
repeated $n_k$ times.   The orbit of $\rho'_\gamma$ has the structure
$Sym^{n_1}\mathbb{C}^{\times}\times \ldots \times Sym^{n_k}\mathbb{C}^{\times}$.
Note that $\beta(\alpha(\rho'_{\gamma})) \in \Omega$.
Consider the orbit of $\rho'_\gamma$.   Each $\rho''$ in the orbit of
$\rho'_{\gamma}$ is such that $\alpha(\rho'')$ is
fixed by $\gamma$, so that $\alpha(\rho'') \in D^{\gamma}$.   
   
Now the centralizer $Z(\gamma) = (\mathbb{Z}/t_1\mathbb{Z}
\wr S_{n_1}) \times \ldots \times (\mathbb{Z}/t_k \mathbb{Z} \wr S_{n_k})$.
But each cyclic group $\mathbb{Z}/t_1\mathbb{Z}, \ldots, \mathbb{Z}/t_k \mathbb{Z}$
acts trivially on $D^{\gamma}$.   Hence $D^{\gamma}/Z(\gamma) = D^{\gamma}/
(S_{n_1} \times \ldots \times S_{n_k}) \cong Sym^{n_1} \mathbb{C}^{\times}
\times \ldots \times Sym^{n_k} \mathbb{C}^{\times}$.

Let $\mathcal{O}(\rho'_{\gamma})$ be the orbit of $\rho'_{\gamma}$.
When we restrict $\alpha$ to this orbit we get
\[
\mathcal{O}(\rho'_{\gamma}) \cong D^{\gamma}/Z(\gamma)
\]
Choose one $\gamma$ in each $\Gamma$-conjugacy class.   We get
\[
\mathcal{G}_\Omega(n) = \bigsqcup \mathcal{O}(\rho'_{\gamma}) =
\bigsqcup D^{\gamma}/Z(\gamma) = \Omega^{+}
\]

If we start with $\rho_1 \otimes sp(r_1) \oplus \ldots \oplus \rho_m \otimes
sp(r_m)$ where $\rho_1, \ldots, \rho_m$ are not all equivalent (after
unramified twist), then we partition this sum in an obvious way and apply the
above argument piecewise. To this end, we note the following  decompositions. 
 Let $X= X_1\times X_2$, $\Gamma =\Gamma_1\times
\Gamma_2$. Then we have 
\[
\widehat{X}/\Gamma = 
 \widehat{X_1}/\Gamma_1 \times \widehat{X_2}/\Gamma_2\; \mbox{\rm and}\; 
\mathcal{A}_{\Omega_1 \times \Omega_2}(n)  = \mathcal{A}_{\Omega_1}(n)  \times
\mathcal{A}_{\Omega_2}(n)
\]
where $\mathcal{A}_{\Omega}(n) = (inf.ch.)^{-1} \Omega$ and $\Omega_1, \Omega_2$ are
disjoint. 

{\sc Example}.  Let $T$ be the diagonal subgroup of $G=GL(2)$ and let
$\Omega$ be the component in $\Omega \, G$ containing
the cuspidal pair
$(T,1)$.     Then $\sigma \in \mathcal{A}(2)$ is {\it arithmetically
unramified} if $inf.ch.\sigma \in \Omega$.   If $\pi_F(\rho') =
\sigma$ then $\rho'$ is a $2$-dimensional representation of
$W'_F$ and there are two possibilities:

{\it $\rho'$ is reducible}, $\rho' = \psi_1 \oplus \psi_2$ with
$\psi_1, \psi_2$ unramified quasicharacters of $W_F$.   So
$\psi_j(w) = z_j^{d(w)}, z_j \in \mathbb{C}^{\times}, j =1,2$.
We have $\pi_F(\rho') = Q(\psi_1, \psi_2)$ where $\psi_1$ does not precede $\psi_2$. 
  In particular we
obtain the 
$1$-dimensional representations of $G$ as follows:
\[
\pi_F(| \;|^{1/2}\psi \oplus | \;|^{-1/2}\psi) = Q(| \;|^{1/2}\psi,
| \;|^{-1/2}\psi) = \psi \circ det
\]
{\it $\rho'$ is indecomposable}, $\rho' = \psi \otimes sp(2)$.
Then $\pi_F(\rho') = Q(\Delta_1)$ with $\Delta_1 = \{\psi, | \;|\psi\}$.
In particular we have $\pi_F(| \;|^{-1/2}\psi \otimes sp(2)) =
Q(\Delta_2)$ with $\Delta_2 = \{|\;|^{-1/2} \psi, | \;|^{1/2} \psi\}$
so $\pi_F(| \;|^{-1/2} \psi \otimes sp(2)) = \psi \otimes St(2)$
where $St(2)$ is the Steinberg representation of
$G$.

The orbit of $(T,1)$ is $D = (\mathbb{C}^{\times })^2$, and $W(T,D) =
\mathbb{Z}/2\mathbb{Z}$.   Then $\Omega \cong
(\mathbb{C}^{\times })^2/\mathbb{Z}/2\mathbb{Z}\cong Sym^2 \,
\mathbb{C}^\times$.   The extended quotient is $\Omega^{+} = Sym^2 \; \mathbb{C}
^{\times} \sqcup \mathbb{C}^{\times}$. 
The twisted projection $\beta$ sends $\{z_1, z_2\}$ to $\{z_1, z_2\}$ and
$z$ to $\{z, q^{-1}z\}$.   Note the twist by $q^{-1}$ where $q$ is the
cardinality of the residue field of $F$.

\mbox{}

\noindent{\bf 2. A geometric model for the Hecke algebra $\mathcal{H}(G)$}

Let $\rho_1, \ldots, \rho_m$ be irreducible representations of $W_F$
such that $det \, \rho_j$ is unitary, $1 \leq j \leq m$, and let
$\rho' =  \rho_1 \otimes sp(r_1) \oplus \ldots \oplus \rho_m \otimes
sp(r_m)$.   Let $\psi_j \in \Psi(W_F)$ and give each $\rho_j$ an unramified
twist as follows:
\[
\rho'' =  \psi_1 | \;\;|^{(1-r_1)/2} \rho_1 \otimes sp(r_1)
\oplus \ldots \oplus \psi_m | \;\;|^{(1 - r_m)/2} \rho_m \otimes sp(r_m)
\]
Then $\rho'$ and $\rho''$ are in the same component of $\mathcal{G}(n)$.

\begin{theorem}\label{3.3}  The map defined by
\[
\bigoplus_j (\psi_j | \;\;|^{(1 - r_j)/2} \rho_j \otimes sp(r_j)) \mapsto
\bigoplus_j (\psi_j |\psi_j|^{-1} | \;\;|^{(1 - r_j)/2} \rho_j \otimes sp(r_j))
\]
determines a deformation retraction of $\mathcal{A}(n)$ onto
the tempered dual $\mathcal{A}^t(n)$.
\end{theorem}
Proof. This depends on the Langlands classification of tempered representations
of $GL(n)$, as in \cite[Prop. 2.2.1]{Ku} and  \cite[p. 384]{Ku}.

\begin{theorem} 
The periodic cyclic homology of $\mathcal{H}(G)$ is isomorphic to the 
periodised de Rham  cohomology  supported on
finitely many  components of the smooth dual of $G$. 
\end{theorem}
Proof. 
We have the Bernstein decomposition for the Hecke algebra $\mathcal{H}(G)$:
$
\mathcal{H}(G) \cong \bigoplus \mathcal{H}(\Omega).
$
We now use \cite[Theorem 7.8]{BHP2} and note that 
de Rham cohomology is invariant under deformation retraction to get
$$
HP_0  (\mathcal{H}(\Omega)) \cong H^{ev}(\Omega^{+}; \mathbb{C})\quad
\mbox{\rm and}\quad
HP_1  (\mathcal{H}(\Omega)) \cong H^{odd}(\Omega^{+}; \mathbb{C})
$$
where the variety $\Omega^+$ is given its classical topology. 
It now follows from Theorem 2 that 
\[
HP_0  (\mathcal{H}(\Omega)) \cong H^{ev}(\mathcal{G}_\Omega(n); \mathbb{C})\quad
\mbox{\rm and} \quad
HP_1  (\mathcal{H}(\Omega)) \cong H^{odd}(\mathcal{G}_\Omega(n); \mathbb{C})
\]

Since $HH_j(\mathcal{H}(\Omega)) = 0$ for all $j>n$, periodic cyclic homology 
commutes with 
direct limits (cf. \cite[Theorem 2]{BP}) and the Theorem follows.

\mbox{}

{\footnotesize
J.B.: School of Mathematical Sciences, University of Exeter, North Park Road, Exeter, EX4 4QE,
U.K., brodzki@maths.ex.ac.uk.

R.J.P.: Department of Mathematics, University of Manchester,
Manchester, M13 9PL, U.K., roger@ma.man.ac.uk.
}
\end{document}